\documentclass[12pt]{amsart}

\newtheorem{theorem}{Theorem}[section]

\newtheorem{defn}{Definition}[section]

\newtheorem{lem}{Lemma}[section]

\newtheorem{cor}{Corollary}[section]

\newtheorem{prop}{Proposition}[section]

\allowdisplaybreaks[4]
\begin{document}

\title{Summability Kernels for $L^p$ Multipliers}


\author{P. Mohanty}


\address{Department of Mathematics, Indian Institute of Technology, 
Kanpur-208016, India}




\email{mparasar@yahoo.com}


\thanks{The first author was supported by CSIR}


\author{S. Madan}

\address{Department of Mathematics, Indian Institute of Technology, 
Kanpur-208016, India}

\email{madan@iitk.ac.in}







\begin{abstract}
In this paper we have characterized the space of summability kernels 
for the case $p=1$ and $p=2$. For other values of $p$ we give a necessary 
condition for a function $\Lambda$ to be a summability kernel.  
For the case $p=1$, we have studied the properties of measures which 
are transferred  from $M(\mathbb Z)$ to $M(\mathbb R)$ by  summability 

kernels. Further,
 we have extended {\bf every} $l_p$ sequence to $L^q(\mathbb R)$ 
multipliers for certain values of $p$ and $q$.
\end{abstract}

\maketitle

\clearpage

\section{Introduction}

In this paper, we extend results proved by Jodeit~\cite{J}, Asmar, 
Berkson, and Gillespie~\cite{ABG1}, Berkson, Paluszy${\rm{\check 
n}}$ski,
 and Weiss~\cite{BPW} in several different ways.

For $1\leq p < {\infty}$, let $M_p(\mathbb R)$ 
(respectively $M_p(\mathbb  Z)$)
 denote  the space of Fourier multipliers for $L^p(\mathbb R)$
 (respectively for $L^p(\mathbb  T)$ ). It is well known that $M_1(\mathbb R)$ can be 

identified with $M(\mathbb R)$, set of all bounded regular measure on $\mathbb R$ and also 

$M_2(\mathbb R)$ can be indentified with $L^\infty(\mathbb R)$. Similar results are true 

for any locally compact abelian group $G$. Here we identify $\mathbb  T$
  with  $[0, 1)$ 
and for $f\in {L^1(\mathbb R)}$ we define its Fourier transform as
$\hat{f}(\xi)=\int\limits_{\mathbb R}f(x)\;e^{-2\pi i\xi x}dx$.

Given a sequence $\phi\in M_p(\mathbb  Z)$, a natural question is when 
and how can $\phi$ be extended to a measurable function $W_\phi$ on 
$\mathbb R$ such that $W_\phi\in M_p(\mathbb R)$. In \cite{J}, 
Jodeit proved that for $1<p<\infty$, the piecewise constant extension 
$\sum\limits_{n\in\mathbb  Z}\phi(n)\chi_{[0,1)}(\xi-n),$ as well as 
the piecewise linear extension do yield multipliers on $L^p(\mathbb 
R)$.
 Fig${\rm\grave a}$-Talamanca and Gaudry~\cite{FG}, using the 
 characterization of the multiplier spaces as dual spaces, proved that 
 the piecewise quadratic extension of $\phi\in M_p(\mathbb  Z)$, 
 $1\leq p<\infty$ is in $M_p(\mathbb R)$. Their method of proof in 
 fact, proves this result for a large class of 
``extensions"~\cite{ABG1}.
  In order to explain our results we need a definition.
\begin{defn}
A measurable function $\Lambda$ on ${\mathbb R}$  is called a 
summability 
kernel for $L^p(\mathbb R)$ multipliers,  if for each 
 $\phi \in M_p(\mathbb  Z)$,  the series 
\begin{equation}
W_{\phi, \Lambda}(\xi)  =  \sum\limits_{n\in\mathbb  Z}
\phi(n)\Lambda(\xi-n)
\label{wphi}
\end{equation}
converges a.e., belongs to $M_p(\mathbb R)$ and there exists a constant 
$C_{p, \Lambda}$ 
such that $\|W_{\phi, \Lambda}\|_{M_p({\mathbb R})} 
\leq C_{p, \Lambda}\|\phi\|_{M_p(\mathbb  Z)}. $
Let $S_p(\mathbb R)$ denote the set of all summability kernels for 
$L^p(\mathbb R)$ multipliers.
\end{defn}

By a different (and powerful) technique, that of transference, Berkson, 
Paluszy${\rm{\check n}}$ski, and Weiss~\cite{BPW} proved that if 
$\Lambda$ is a measurable function on $\mathbb R$ with compact 
support then $\Lambda$ is a summability kernel if and only if 
$\Lambda\in M_p(\mathbb R)$. Note that if $supp~\Lambda$ is compact, 
the series in (\ref{wphi}) converges for every $\phi\in M_p(\mathbb  
Z)$.
 Without such a hypothesis, it is not clear that (\ref{wphi}) 
converges.
  We begin with a lemma which partially explains the a.e. convergence 
of
    (\ref{wphi}). Consider the following set:\\
$S_p^0(\mathbb R)=\Bigl\{\Lambda\in L^\infty(\mathbb R)$:\;For  
each finitely supported $\phi \in M_p(\mathbb  Z),$ the  
function
$W_{\phi,\Lambda}(\xi)=\sum\limits_{n\in\mathbb  Z}\phi(n)
\Lambda(\xi-n)$
belongs to $M_p(\mathbb R)$ and there 
exists a
 constant $C_{p, \Lambda}$ such that
 $\|W_{\phi, \Lambda}\|_{M_p(\mathbb R)} \leq C_{p, \Lambda}\|
\phi\|_{M_p(\mathbb  Z)}\Bigr\}$.

\begin{lem} 

(i)For $1<p<\infty$, let $\Lambda \in S_p^0({\mathbb R})$. If for
$\phi\in M_p(\mathbb  Z)$, the series
$$\sum\limits_{n\in\mathbb  Z}\phi(n)\Lambda(\xi-n) = W_{\phi, \Lambda}
(\xi)$$
converges a.e.,  then $\Lambda\in S_p(\mathbb R).$\\
(ii) If $p=1,\; \Lambda \in S_1^0(\mathbb R)$ and
$\delta_{\Lambda}=\sup\limits_{\xi} \sum\limits_{n\in\mathbb  Z}|
\Lambda(\xi+n)|<\infty$
then $\Lambda \in S_1(\mathbb R)$.

\label{lem:1.2}
\end{lem}
\noindent {\bf Proof}
  (i)For a.e. $\xi$, we have
$$W_{\phi, \Lambda}(\xi)=\lim_{N\rightarrow \infty}\ \sum\limits_
{-N}^N\phi(n)\Lambda(\xi-n)=\lim_{N\rightarrow \infty}\ W_{\phi_{N},
 \Lambda}(\xi)$$
where $\phi_N(n) = \left\{\begin{array}{cl}

                            \phi(n) & \mbox{if $|n| \leq N$} \\

                            0       & \mbox{otherwise}.

                         \end{array}\right. $
\noindent 

Since $1<p<\infty$,  $\phi_N\in M_p(\mathbb  Z), $ and 
$\|\phi_{N}\|_{M_p(\mathbb  Z)}\leq C_p\|\phi\|_{M_p(\mathbb  Z)}$,  
where $C_p$ is a constant independent of $N$. This follows from 
 M.Riesz theorem ~\cite{EG}.  Since $\Lambda \in S_p^0(\mathbb 
R)$ we have
\begin{eqnarray*}
\|W_{\phi_{N}, \Lambda}\|_{M_p(\mathbb R)}
&\leq&   C_{p, \Lambda} \|\phi_{N}\|_{M_p(\mathbb  Z)}\\
 & \leq & C_p C_{p, \Lambda}\|\phi\|_{M_p(\mathbb  Z)}.
\end{eqnarray*}
Hence $W_{\phi_{N}, \Lambda}(\xi)\rightarrow W_{\phi, \Lambda}(\xi)$ 
pointwise a.e. and boundedly , so $W_{\phi, \Lambda}\in 
M_p(\mathbb R)$.

 (ii) For $p=1$,  the additional condition guarantees the convergence 
of $\sum\limits_{n\in\mathbb  Z}\phi(n)\Lambda (\xi~-~n~)$ for every 
$\phi\in M_1(\mathbb  Z)$. If $K_N$ is the $N$th F$\acute{\mbox e}$jer
 kernel on $\mathbb R$,
  let $\phi_N=\hat{K}_N\phi$. Then again 
$W_{\phi_{N}, \Lambda}(\xi)\rightarrow  W_{\phi, \Lambda}(\xi)$ 
pointwise a.e. and boundedly. Hence, $\Lambda\in S_1(\mathbb R)$.

Clearly $S_p(\mathbb R)\subseteq S_p^0(\mathbb R)$. In \S 2 of this 
paper we  characterize the space $S_p^0(\mathbb R)$ for $p=1$ and for 
$p=2$. Further we  prove fairly general results which make precise the 
reasons why summability kernels allow the transference of multipliers from 
$L^p(\mathbb  T)$ to $L^p(\mathbb R)$.

In \S 3 we restrict ourselves to the case $p=1$. We investigate some 
properties of measures which are stable under transference by summability 
kernels, such as the properties of being  discrete, continuous or 
absolutely continuous.

In \S 4 we will show that for some values of $p$ and $q$ every sequence 
in $l_p(\mathbb  Z)$ can be extended  to an $L^q(\mathbb R)$ multiplier 
by means of suitable summability kernels.

\section{Summability Kernels for $L^p({\mathbb R}),  1\leq 
p<\infty$}\label{sec:2}

For a function $\Lambda\in L^\infty(\mathbb R)$ denote 
$\delta_\Lambda={\rm{ess}}\sup\limits_\xi\sum\limits_{n\in\mathbb  Z}|\Lambda(\xi+n)|$ 
.  

 In the following theorem we  characterize $S_p^0(\mathbb R)$, for 
$p=1$ and $p=2$.

\begin{theorem}
(i) $S_2^0=S_2=\{\Lambda\in L^\infty(\hat{\mathbb 
R}):\delta_\Lambda<\infty\}$.\\
(ii) $S_1^0=\{\Lambda\in L^1(\mathbb R)^\wedge:\Lambda=\hat 
F\;{\rm{with}}\;\delta_F<\infty\}.$

\label{thm:0.2}

\end{theorem}
\noindent{\bf Proof:}

 (i) Let $\Lambda\in L^\infty(\mathbb R)\;{\rm{and\;suppose}}\;\delta_
 \Lambda<\infty$. Then clearly $\Lambda\in S_2(\mathbb R)$. 
 Let $\Lambda\in S_2^0(\mathbb R)$. For $\xi\in\mathbb R$, let 

$\phi_\xi(n)=e^{i\theta}$ where $\Lambda(\xi-n)=e^{i\theta}|\Lambda(\xi-n)| \;\;(\theta$ of 

course depends on $\xi$ and $n)$.

 Then $\|\phi_\xi\|_{M_2(\mathbb Z)}=1$. Also for every positive integer $N$ we can find a 

null subset $E_N$ of $\mathbb R$ such that 

$$\sum\limits_{|n|\leq N} \phi_\xi(n)\Lambda(\xi-n)\leq \|
\sum\limits_{|n|\leq N}\phi_.(n)\Lambda(.-n)\|_\infty\;\;\;\;
\forall \xi\in \mathbb R\setminus E_N.$$

For $\xi\in \mathbb R\setminus\cup_N E_N$, we have 

\begin{eqnarray*}
\sum\limits_{n\in\mathbb  Z}|\Lambda(\xi-n)|& = & 
\sup\limits_N\sum\limits_{|n|\leq N}|\Lambda(\xi-n)|\\
& = & \sup\limits_N\sum\limits_{|n|\leq 
N}\phi_{\xi}(n)\Lambda(\xi-n)\\
& \leq & \sup\limits_N\|\sum\limits_{|n|\leq 
N}\phi_.(n)\Lambda(.-n)\|_{\infty}\\
& \leq & C_{2,\Lambda}<\infty\hspace{0.75in}({\rm{as}}\;\Lambda\in S_2^0(\mathbb R)).
\end{eqnarray*}
So $\delta_\Lambda<\infty$. We already know that 
$S_2(\mathbb R)\subseteq S_2^0(\mathbb R)$. Hence 
$S_2(\mathbb R)=S_2^0(\mathbb R)$.\\
(ii) Let $\Lambda=\hat F$, where $F\in L^1(\mathbb R)$ and 
$\delta_F<\infty$. For a finite sequence $\{\phi(n)\}$,  let 
$P(x)=\sum\limits_n\phi(n)e^{2\pi inx}$. Then 

\begin{eqnarray*}
W_{\phi, \Lambda}(\xi) & = & \sum\limits_{n\in\mathbb  Z}\phi(n)\hat 
F(\xi-n)\\
& = & \sum\limits_n\phi(n)\int_\mathbb R F(x)e^{-2\pi ix(\xi-n)} dx\\
& = & (FP^\#)^\wedge(\xi),
\end{eqnarray*}

where $P^\#$ is the 1-periodic extension of $P$ on $\mathbb R$. 

 We have $FP^\#\in L^1(\mathbb R)$ as $P^\#\in L^\infty(\mathbb R)$. 
 Thus 

\begin{eqnarray*}
\|W_{\phi, \Lambda}\|_{M_1(\mathbb R)} & =  & \|FP^\#\|_1\\
& = & \int_0^1\sum\limits_n|F(x+n)||P(x)| dx\\
& \leq & \delta_F\|P\|_{L^1(\mathbb  T)}=\delta_F\|\phi\|_{M_1(\mathbb  
Z)}.
\end{eqnarray*}

Hence $\Lambda\in S_1^0(\mathbb R)$.

Conversely, suppose $\Lambda\in S_1^0(\mathbb R)$. Then taking 
$\phi(n)=\delta_{n, 0}$ we have $\Lambda\in M_1(\mathbb R)=
M(\mathbb R)^\wedge$.  So $\Lambda=\hat\mu$ for some 
$\mu\in M(\mathbb R)$. For a finite sequence $\{\phi(n)\}$ and 
$P$ as above we have 

$$W_{\phi, \Lambda}(\xi)=(P^\#\mu)^\wedge(\xi),$$

where $P^\#\mu$ denotes the measure given by 
$d(P^\#\mu)(x)=P^\#(x)d\mu(x).$ Then 

$$\|P^\#\mu\|_{M(\mathbb R)}  =  \|W_{\phi, \Lambda}\|_{M_1(\mathbb R)}
\leq  C_{1, \Lambda}\|\phi\|_{M_1(\mathbb  Z)}
 =C_{1, \Lambda}\|P\|_{L^1(\mathbb  T)}.$$ For $k\in\mathbb  Z$ let 
$\mu_k=\tau_{-k}(\mu|_{[k, k+1)})$, i.e., $\mu_k$ is a measure 
supported on $[0,1)$. We have  
$\|P\mu_k\|_{M(\mathbb  T)}\leq C_{1, \Lambda}\|P\|_{L^1(\mathbb  T)}$. 
By density of trigonometric polynomials, $\|f\mu_k\|_{M(\mathbb  T)}\leq C_{1,\Lambda}\|f\|_{L^1
(\mathbb  T)}$ for $f\in L^1(\mathbb T)$. It follows 
easily that each $\mu_k$ is absolutely continuous, hence so is $\mu$. 
Let $d\mu=F dx$ with $F\in L^1(\mathbb R)$, so that  $\Lambda=\hat F$. 
Now 

\begin{eqnarray*}
\|W_{\phi,\Lambda}\|_{M_1(\mathbb R)}=\|P^\#F\|_{L^1(\mathbb R)} & = & 
\|P|F|^\#\|_{L^1(\mathbb  T)}\\
& \leq & C_{1,\Lambda}\|\phi\|_{M_1(\mathbb  
Z)}\\
& = & C_{1,\Lambda}\|P\|_1.
\end{eqnarray*}

Therefore $|F|^\#$ defines a continuous linear functional on 
$L^1(\mathbb  T)$, so by duality
$|F|^\#\in L^\infty(\mathbb  T)$, and 
$\||F|^\#\|_{ L^\infty(\mathbb  T)}=\delta_F<\infty$.

\noindent{\bf{Remarks:}}

The two conditions appearing in Theorem~\ref{thm:0.2} for $p=1$ 
and for $p=2$ seem to be very different. We will analyse these further 
to obtain a more unified formulation. If $p=1$,  the condition $\delta_F<\infty$

 is equivalent to 
saying that for a.e. $x\in\mathbb R$, the sequences 
$\{F(x+n)\}_{n\in\mathbb  Z}$ define, by convolution, operators on 
$l_1(\mathbb  Z)$ or in other words, we have \\
$\rm(a_1)$ $\Lambda_x^\#\in M_1(\mathbb  T)$ 
( where $(F(x+.))^\wedge=\Lambda_x^\#$).

If $p=2$. The condition on $\Lambda$, namely $\delta_\Lambda<\infty$ 
 implies that\\
${\rm(a_2)}$ for a.e. $x,\;\Lambda_x^\#\in L^\infty(\mathbb  T)=
M_2(\mathbb  T)$, where $\Lambda_x=e^{2\pi ix.}\Lambda.$

From  $\rm(a_1)$ and $\rm(a_2)$ above we get a  condition which we  
show
 is necessary for $\Lambda$ to be  a summability kernel for 
 $L^p(\mathbb R)$ multipliers. 
For $1\leq p\leq 2$ define\\

${\mathcal{F}}_p=\{\Lambda\in L^\infty(\mathbb R): 
{\rm{for\;\;a.e.}}\;x\in[0, 1),\;\;\Lambda_x^\#\in M_p(\mathbb  T)\;
{\rm{with}}\;\|\Lambda_x^\#\|_{M_p(\mathbb  T)}\in L^\infty[0, 1)\}$.
\begin{prop}
$S_p(\mathbb R)\subseteq {\mathcal{F}}_p\;\;{\rm{for}}\;1\leq p\leq 2.$ 
\label{prop:2.2}
\end{prop}
\noindent{\bf{Proof:}}

 Let $\Lambda\in S_p(\mathbb R)$, then for every $x\in\mathbb R, \Lambda_x\in S_p(\mathbb R)$. 
 In fact 
\begin{eqnarray*}
 W_{\phi,\Lambda_x}(\xi) & = & \sum\limits_n \phi(n)\Lambda_x(\xi-n)\\
& = & e^{2\pi ix\xi}\sum\limits_n\phi(n)e^{-2\pi ixn}\Lambda(\xi-n).
\end{eqnarray*}
Since  $\phi_x(n)=e^{-2\pi ixn}\phi(n)$ belongs to $M_p(\mathbb  Z)$ 
whenever $\phi\in M_p(\mathbb  Z)$ with equal norm, we have 
$ W_{\phi,\Lambda_x}\in M_p(\mathbb R)$ and 
\begin{eqnarray*}
\| W_{\phi,\Lambda_x}\|_{M_p(\mathbb R)} & = &\| W_{\phi_x,\Lambda}\
|_{M_p(\mathbb R)}\\
& \leq & C_{p, \Lambda}\|\phi\|_{M_p(\mathbb  Z)}.
\end{eqnarray*}
 Now take $\phi(n)\equiv 1$. Then $ W_{1,\Lambda_x}(\xi)=\sum
 \limits_n e^{-2\pi ix(\xi-n)}\Lambda(\xi-n)$ belongs to 
 $M_p(\mathbb R)$ and is 1-periodic.  Thus by de Leeuw's 
 result~\cite{DL}, $W_{1, \Lambda_x}\in M_p(\mathbb  T)$ and 
 $\|W_{1, \Lambda_x}\|_{M_p(\mathbb  T)}\leq C_{p, \Lambda}\;
 {\rm{for\;all}}\; x\in[0, 1)$. Therefore 
 $\|{\Lambda_x}^\#\|_{M_p(\mathbb  T)}\in L^\infty [0, 1)$.

 In ~\cite{ABG1}, Asmar, Berkson, and Gillespie by using transference technique and considering 

the $M_p(\mathbb R)$ as the dual of 

Figa-Talamanca and Herz algebra $A_p(\mathbb R)$ proved that if $J\in L^1(\mathbb R)$ with 

$\hat J$ has compact support, $J$ is absolutely continuous and $J^\prime\in L^2(\mathbb R)$ 

then $\hat J$ and $J$ belong to $S_p(\mathbb R)$. By suitably refining 
the proof (Theorem III.4 in ~\cite{ABG1}), we have the following theorem, 
which  improves on their result.  
\begin{theorem}
(i) If $\Lambda_1,\;\Lambda_2\in{\mathcal{F}}_p\cap A(\mathbb R)$ then 
$\Lambda=\Lambda_1\Lambda_2\in S_p^0(\mathbb R).$\\

 (ii) If, in addition, either $\delta_{\Lambda_1}<\infty$ or 
 $\delta_{\Lambda_2}<\infty$ then $\Lambda\in S_p(\mathbb R)$.
\label{thm:1.2}
\end{theorem}

The proof is exactly as in ~\cite{ABG1}, except for the improvement in their inequallity 

III.12, where we use the better estimate from the transference principle.  For details see 

~\cite{PM}. In particular $S_1^0.S_1^0\subset S_p$.

It is clear that there is a large class of 
summability kernels which do not have compact support nor are Fourier 
transforms of  compactly supported integrable functions.
Let $\Lambda\in M_p(\mathbb R)$ and $\Lambda_k=\chi_{[k,k+1)}\Lambda$. If 

 $\sum\limits_k\|\Lambda_k\|_{M_p(\mathbb R)}<\infty$ then 

 $W_{\phi, \Lambda}\in M_p(\mathbb R)$. It follows that every such class of function 

is a summability kernel. 
  \section{Continuous, Absolutely continuous, and Discrete measures}
For the case $p=1$,  $M_1(\mathbb R)$ and $M_1(\mathbb  Z)$ are 
identified with the set of bounded regular measures 
$M(\mathbb R)$ and $M(\mathbb  T)$ respectively.   So if 
$\phi\in M_1(\mathbb  Z)$ then 
$\phi=\hat{\nu}$ for some $\nu\in M(\mathbb  T)$, and if 
$\Lambda$ is a
 summability kernel we have $W_{\phi, \Lambda}~\in M_1(\mathbb R)$. 
 i.e. $W_{\phi, \Lambda}=\hat{\mu}$ for some $\mu\in M(\mathbb R)$. 

Here  we study some properties of measures which are carried over from $\nu$ to 
$\mu$.
  Let
${\mathcal F}_0=\{\Lambda\in S_1^0(\mathbb R):\delta_\Lambda<\infty\}$.

\begin{theorem}

Let $\Lambda \in {\mathcal F}_0, \ \nu\in M({\mathbb  T})$ and define
$\hat{\mu}(\xi)=W_{\hat{\nu},
\Lambda}(\xi)=\sum\limits_{n}\hat{\nu}(n)\Lambda(\xi-n)$,  (here
$\mu\in M(\mathbb R)$).   \\
 (a) If $\nu$ is an absolutely continuous measure on ${\mathbb  T}$, 
 then $\mu$ is an absolutely continuous measure on $\mathbb R$
 (both with respect to the Lebesgue measure).\\
(b) If $\nu$ is a discrete measure,  then either $\mu\equiv 0$ or $\mu$ 
is a discrete measure.

\label{thm:3.2}

\end{theorem}

\noindent {\bf Proof:}

 (a) 
First assume that $d\nu(x)=P(x)dx$,  where $P$ is a trigonometric 
polynomial.   Then
\begin{eqnarray*}
\hat{\mu}(\xi) & = & \sum\limits_{n}\hat{P}(n)\Lambda(\xi-n)\\
& = & \sum\limits_{n}\hat{P}(n)\ \left(e^{ 2\pi 
in.}g\right)^{\wedge}(\xi)\\
& = & (P^\#g)^\wedge(\xi),
\end{eqnarray*}
where $\Lambda = \hat{g},\;\;g\in L^1(\mathbb R)$. Let $h=P^\# g$. 
Then, 
\begin{equation}
 \| h \|_{L^1(\mathbb R)} = \|\hat{\mu}\|_{M_1(\mathbb R)}\leq 
C_{1,\Lambda}\|\hat{\nu}\|_{M_1(\mathbb  Z)}=C_{1,\Lambda}\|P\|_{L^1(\mathbb  
T)}.
\label{eqn:6.2}
\end{equation}
Now if $\nu$ is an absolutely continuous measure on $\mathbb  T$,  let 
$\hat{\nu}=\hat{F}$ for $F\in L^1(\mathbb  T)$.   There exists a 
sequence 
$\{P_N\}$ of trigonometric polynomials such that $P_N\rightarrow F$ in 
$L^1(\mathbb  T)$.  For each $P_N$ define $h_N$ as above.
 Then from Eqn.~(\ref{eqn:6.2})
$$\|h_N-h_M\|_{L^1(\mathbb R)}\leq C_{1, 
\Lambda}\|P_N-P_M\|_{L^1(\mathbb  T)}. $$
Let $h_N\rightarrow h$ in $L^1(\mathbb R)$.   Now
\begin{eqnarray*}
|\hat{h}_N(\xi)-\hat{\mu}(\xi)| & \leq & 
\sum\limits_{n}|\hat{P}_N(n)-\hat{F}(n)|\; |\Lambda (\xi-n)|\\
& \leq & \|P_N-F\|_{L^1(\mathbb R)}\ \delta_{\Lambda}.
\end{eqnarray*}
So, $\hat h=\hat \mu$. Hence, $d\mu(x)=h(x)dx.$

(b) Let $\nu=\sum\limits_{j=1}^{\infty}\ \alpha_j\delta_{x_j}$ 
be a discrete
measure on ${\mathbb  T}$,  $x_j\in {\mathbb  T}$ and 
$\sum\limits_{j}|\alpha_j|<\infty$. 
  If $\Lambda=\hat{g}$,  with $g\in L^1(\mathbb R)$, ( by Fourier inversion 
we may assume that $g$ is continuous) and  if $\tilde g(x)=g(-x)$,  then
\begin{eqnarray*}
\hat\mu(\xi) & = & \sum\limits_{n}\ \hat{\nu}(n)\Lambda(\xi-n)\\
& = & \sum\limits_{n}\ \left(\sum\limits_j\alpha_j e^{2\pi 
inx_j}\right)\ \left(e^{2\pi
i\xi\cdot}\ \tilde{g}\right)^{\wedge}(n)\\
& = & \sum\limits_j\alpha_j\sum\limits_{n\in \mathbb  
Z}\left(\tau_{-x_j}\ e^{2\pi i\xi\cdot}\tilde{g}\right)(n) \\
& = & \sum\limits_j\ \alpha_j\sum\limits_{n\in \mathbb  Z}\ e^{2\pi 
i\xi(n+x_j)}\ \tilde g(x_j+n).
\end{eqnarray*}
The last but one equality follows 
by Poisson summation formula, since $\Lambda\in S_1^0(\mathbb R)$ will 
imply $\sum\limits_n|g(x+n)|<\infty$ and $\Lambda\in {\mathcal F}_0$ 
will imply $\sum\limits_{n\in\mathbb  Z}|\hat g(\xi+n)|<\infty$.    
Now it is clear that
$$\mu=\sum\limits_{n\in \mathbb  Z}\ \sum\limits_j\alpha_j\tilde 
g(x_j+n)\ \delta_{x_{j}+n}.$$
Hence $\mu$ is either the zero measure or is  discrete.

To consider similar results for continuous measures,  we need some 
additional conditions on the summability kernel,  and this is the 
content
 of the following two results, both of which use Wiener's 
lemma~\cite{K}.
\begin{theorem}
Let $\Lambda \in S_1(\mathbb R)$,  and suppose that $supp\; \Lambda$ is 
compact.   Then if $\nu$ is a continuous measure,  so is $\mu$.
\label{thm:4.2}
\end{theorem}
\noindent {\bf Proof:}
By Wiener's Lemma 
\begin{eqnarray*}
\mu\{y\} & = & \lim_{\lambda\rightarrow\infty}\ \frac{1}{2\lambda}\ 
\int_{-\lambda}^{\lambda}\ \hat{\mu}(\xi)\ e^{2\pi i\xi y}\ d\xi\\
& = & \lim_{\lambda\rightarrow\infty}\ \frac{1}{2\lambda}\ 
\sum\limits_{n\in \mathbb  Z}\hat{\nu}(n)\ \int_{-\lambda}^{\lambda}\ 
\Lambda(\xi-n)\ e^{2\pi i\xi y}\ d\xi\\
& = & \lim_{\lambda\rightarrow\infty}\ \frac{1}{2\lambda}\ 
\sum\limits_{n}\hat{\nu}(n)\ I_{\lambda}^n(y),
\end{eqnarray*}
where $$I_{\lambda}^n(y) = \int_{-\lambda}^{\lambda}\ \Lambda(\xi-n)\ 
e^{2\pi i \xi y}\ d\xi. $$
Let $supp\; \Lambda \subseteq [-N, N]. $  Then for each $\lambda >0$,  
if $|n|> N+\lambda$,  clearly 
$I_{\lambda}^n(y)=0\;\forall y\in \mathbb R$.   Now let $\lambda > 2N.$

\noindent {\bf Case 1}
Suppose $|n|\leq \lambda -N$,  then
\begin{eqnarray*}
I_{\lambda}^n(y) & = & \int_{-\lambda+n}^{\lambda+n}\ \Lambda (\xi)\ 
e^{2\pi i (\xi+n)y}\ d\xi\\
& = & \hat{\Lambda}(-y)\ e^{2\pi i yn}
\end{eqnarray*}
since $[-N, N]\subset [-\lambda+n, \lambda+n]. $

\noindent{\bf Case 2}
$\lambda-N\leq |n|\leq \lambda+N$.   Then
$$I_{\lambda}^n(y)=\int_{-N+n}^{\lambda}\ \Lambda(\xi-n)\ e^{2\pi i \xi 
y} \ d\xi.$$
So, in both the cases 
$$|I_{\lambda}^n (y)|\leq \|\Lambda\|_{L^1(\mathbb R)}. $$
Hence,  for $\lambda > 2N$
\begin{eqnarray*}
\mu\{y\} & = & \lim_{\lambda\rightarrow \infty}\ \frac{1}{2\lambda}\ 
\left(\sum\limits_{|n|\leq \lambda-N}\ \hat{\nu}(n)\ 
I_{\lambda}^n(y)+\sum\limits_{N+\lambda\geq |n|\geq \lambda-N}\ 
\hat{\nu}(n)I_{\lambda}^n(y)\right)\\
& = & \lim_{\lambda\rightarrow \infty}\ \frac{1}{2\lambda}\ 
\sum\limits_{|n|\leq [\lambda-N]}\ \hat{\nu}(n)\ e^{2\pi iyn}\ \hat{\Lambda}(-y) 
\\
&  & + \lim_{\lambda\rightarrow \infty}\ \frac{1}{2\lambda}\ 
\sum\limits_{N+\lambda\geq |n|>\lambda-N}\ \hat{\Lambda}(n)\ I_{\lambda}^n(y).
\end{eqnarray*}
The second limit is zero since the terms are bounded and the number of 
terms is at most $2N$. Applying Wiener's lemma for the continuous 
measure $\nu$ on $\mathbb  T$ for the first limit we have 
$$\mu\{y\} =  \nu\{y_0\}\hat{\Lambda}(-y)=0 \;{\mbox{ where 
}}\;{y_0}\in[0, 1)  \;s. t. \; y=y_{0}+2\pi l  
\;{\mbox { for some }}\;l\in{\mathbb  Z}.$$
The hypothesis that $supp\; \Lambda$ be compact may be too restrictive.  
 It can be replaced by the existence of a suitable decreasing radial 
$L^1$ - majorant $\Lambda_1$, i. e., a function $\Lambda_1$ 
satisfying\\
(a) $\Lambda_1$ is decreasing and radial\\
(b) $\Lambda_1\in L^1(\mathbb R)$\\
(c) $|\Lambda(\xi)|\leq |\Lambda_1(|\xi|)|$.\\

Note that for decreasing, radial $L^1$- function $\Lambda$, 

$\sup\limits_\xi\sum\limits_n |\Lambda(\xi+n)|<\infty$.
\begin{theorem}
Suppose $\Lambda\in {\mathcal F}_0$ and that $\Lambda$ has a decreasing
 radial $L^1$ - majorant $\Lambda_1$.   Then $\mu$ is a continuous 
 measure if $\nu$ is.
\label{thm:5.2}
\end{theorem}
\noindent {\bf Proof:}
Once again,  we use Wiener's lemma.\\
Let
\begin{eqnarray*}
I_{\lambda} & = & \frac{1}{2\lambda}\ \int_{-\lambda}^{\lambda}\ 
|\hat{\mu}(\xi)|^2\ d\xi\\
& = & \frac{1}{2\lambda}\ \int_{-\lambda}^{\lambda}\ 
|\sum\limits_{n}\hat{\nu}(n)\Lambda(\xi-n)|^2\ d\xi\\
& \leq  & \frac{1}{2\lambda}\ \int_{-\lambda}^{\lambda}\ 
\left(\sum\limits_{n}|\hat{\nu}(n)|^2|\Lambda(\xi-n)|\right)\ 
\left(\sum\limits_{n}|\Lambda(\xi-n)|\right)\ d\xi\\
& \leq & \delta_{\Lambda}\left[\frac{1}{2\lambda}\ \sum\limits_{|n|\leq 
2\lambda}\ |\hat{\nu}(n)|^2\ \int_{-\lambda}^{\lambda}|\Lambda(\xi-n)| 
d\xi \right. \\
& & \qquad + \left. \frac{1}{2\lambda}\ \sum\limits_{|n|>2\lambda}\ 
|\hat{\nu}(n)|^2
\int_{-\lambda}^{\lambda}|\Lambda(\xi-n)| d\xi \right]\\
& = & \delta_{\Lambda}(I_1+I_2), \; \mbox{say}.
\end{eqnarray*}
Now,
\begin{eqnarray*}
I_2 & \leq & \frac{1}{2\lambda}\ 
\sum\limits_{|n|>2\lambda}|\hat{\nu}(n)|^2\ \int_{-\lambda}^{\lambda}\Lambda_1(|\xi-n|)\ d\xi\\
& \leq & \frac{1}{2\lambda}\ \|\hat{\nu}\|_{\infty}\ 
\left(\sum\limits_{n>2\lambda}\
\Lambda_1(|\lambda-n|)+\sum\limits_{n<-2\lambda}\ 
\Lambda_1(|\lambda-n|)\right)\\
&  \rightarrow & 0\;\;\;\;{\rm{as}}\;\lambda\rightarrow\infty,
\end{eqnarray*}
Hence,  using Wiener's lemma for 
$\mathbb  T$ we get
\begin{eqnarray*}
\lim_{\lambda\rightarrow\infty}\ I_{\lambda} & \leq & \delta_{\Lambda}\ 
\lim_{\lambda\rightarrow \infty}\frac{1}{2\lambda}\ 
\sum\limits_{|n|\leq [2\lambda]}\ |\hat{\nu}(n)|^2\ 
\int_{-\lambda}^{\lambda}|\Lambda(\xi-n)|\ d\xi\\
& \leq & \delta_{\Lambda}\|\Lambda\|_{L^1(\hat {\mathbb R})}\ 
\lim_{\lambda\rightarrow\infty}\frac{1}{2\lambda} \sum\limits_{|n|\leq 
[2\lambda]}|\hat{\nu}(n)|^2\\
& = & 0.
\end{eqnarray*}
\section{Extensions of Sequences to $L^p$ Multipliers}
In this section we study a different kind of extension. In  the 
existing literature  the emphasis has been on the  extensions of 
$L^p(\mathbb  T)$-multipliers to $L^p(\mathbb R)$-multipliers. Here we 
will show that for some values of $p$ and $q$  every sequence in 
$l_p(\mathbb  Z)$ can be extended to an $L^q(\mathbb R)$ multiplier by 
means of suitable summability kernels. The idea of our  extension comes 
from the following result of Jodiet ~\cite{J}.
\begin{theorem}  
Let $S\in L^1(\mathbb R)$, $supp\; S\subseteq 
[\frac{1}{4},\frac{3}{4}]$  and suppose its 
 1-periodic extension $S^\#$ from $[0,1)$  has an absolutely summable 
Fourier series. Then
\begin{eqnarray}
W_{\phi,\hat S}(\xi) = \sum\limits_{n\in\mathbb  Z}\phi(n)\hat S(\xi-n)
\end{eqnarray}
is in $M_p(\mathbb R)$ whenever $\phi\in M_p(\mathbb  Z)$ and its norm 
is bounded by 
$C_p\tau \|\phi\|_{M_p(\mathbb  Z)}$ where 
$\tau = \sum\limits_n|(S^\#)^\wedge(n)|$
 and $C_p$ is a constant which depends only on $p$.
\label{thm:one3}
\end{theorem}
It is natural to ask what happens if we assume 
$(S^\#)^\wedge\in l_p(\mathbb  Z)$ for 
$1<p<\infty$. In this case, it follows from Lemma \ref{lem:one3} below  
and H${\rm{\ddot o}}$lder's inequality that the above sum converges for 
every sequence $\{\phi(n)\}\in l_{p^\prime}(\mathbb  Z)$ where 
$\frac{1}{p}+\frac{1}{p^\prime}=1$, and defines a function 
$W_{\phi, \hat S}$ in $L^\infty(\mathbb R).$
\begin{theorem}
 Let $S\in L^1(\mathbb R)$,  
 $ supp\;S\subseteq [\frac{1}{4},\frac{3}{4}]$  and
$\sum\limits_n|(S^\#)^\wedge(n)|^p <\infty$ for $1<p<\infty$. Define 
$W_{\phi, \hat S}(\xi)=\sum\limits_{n\in \mathbb  Z}
{\phi(n)\hat S(\xi-n)}$ for $\phi\in l_{p^\prime}(\mathbb  Z)$. Then
\begin{eqnarray*} 
W_{\phi, \hat S}\in {M_q(\mathbb R)}\;\; for 
\begin{cases}
q\in [\frac{2p}{3p-2},\frac{2p}{2-p}]\;\; if \;1<p<2 \\ 
q\in [\frac{2p}{p+2},\frac{2p}{p-2}]\;\; if\; 2<p<\infty .
\end{cases}
\end{eqnarray*}
 For $p=2,W_{\phi, \hat S}\in {M_q(\mathbb R)}$
 for all $1\leq q <\infty$. Moreover, 
 $$\|W_{\phi, \hat S}\|_{M_q(\mathbb R)}\leq C\tau_p\|\phi\|_p$$
where $\tau_p=(\sum\limits_n|(S^\#)^\wedge(n)|^p)^{\frac{1}{p}}$ and 
$C$ is a constant which depends only on $p$.
\label{thm:two3}
\end{theorem}
By putting a further restriction on $S$ we will get 
$W_{\phi,\hat S}\in M_q(\mathbb R)$ for 
$1\leq q<\infty$ whenever $\phi\in l_p(\mathbb  Z)$, for $1<p<2$. 

For the proof of Theorem \ref{thm:two3}, we first prove a lemma.
\begin{lem}
Let $S\in L^1(\mathbb R)$, $supp\; S\subseteq[\frac{1}{4},\frac{3}{4}]$ 
 and 
$\sum\limits_n|(S^\#)^\wedge(n)|^p<\infty$ for $1\leq p<\infty$. Then 
$\sup\limits_\xi\sum\limits_n|\hat S(\xi+n)|^p < \infty$.\\
Moreover for $p=1$, $\sum\limits_{n\in\mathbb  Z}\hat S(\xi+n)=C$ for 
all $x\in\mathbb R$, where $C$ is a constant.
\label{lem:one3}
\end{lem}
\noindent{\bf Proof:}
Let $\rho\in {C_c}^\infty(\mathbb R)$ be such that $\rho(x)=1$ on 
$[\frac{1}{4},\frac{3}{4}]$,
 $supp\;\rho\subseteq [\frac{1}{8},\frac{7}{8}]$ and 
$|\rho(x)|\leq ~~1\;\;\forall x\in\mathbb R$. For $\xi\in [0,1]$ define
$h_\xi(x) = e^{-2\pi i\xi x}\rho(x)$. Then 
$h_\xi\in {C_c}^\infty(\mathbb R)$  and 
$$({h_\xi}^{(2)})^\wedge(y)  = (2\pi iy)^{2}\hat h_\xi(y)$$
Hence
\begin{eqnarray*}
 |\hat{h_\xi}(y)|& = &\frac{| (h_\xi^{(2)})^\wedge(y)|}{|2\pi 
y|^2}\hspace{0.3in}{\rm{(for}}\;y\not = 0)\\
& \leq &
 \frac{\|h_\xi^{(2)}\|_1}{|2\pi y|^2}\leq\frac{C}{|y|^2}
\end{eqnarray*} 
where the constant $C$ is 
 independent of $\xi$. So in particular  
 $|\hat h_\xi(n)|\leq \frac{C}{n^2}$ for $n\not =0$. Now define 
 $g_\xi(x) = h_\xi(x)S(x)$. Then 
$supp\; g_\xi\subseteq [\frac{1}{4},\frac{3}{4}]$ and
$$\hat{g^\#_\xi}(n)=\hat 
g_\xi(n)=\int_{\frac{1}{4}}^{\frac{3}{4}}e^{2\pi i \xi x}S(x)e^{-2\pi i\xi n} = \hat S(\xi+n)$$
where $g_\xi^\#$ is the 1-periodic extension of $g_\xi$ 
given by $\sum\limits_{n\in\mathbb  Z} g_\xi(x+n)$. Also
 $\hat{g_\xi^\#}(n)=\hat{h_\xi^\#}*(S^\#)^\wedge(n)$. Since  
$\hat{h_\xi^\#}\in l_1(\mathbb  Z)$ and 
$({S^\#})^\wedge\in l_p(\mathbb  Z)$, it follows that 
$\hat {g_\xi^\#} \in l_p(\mathbb  Z)$ for $1\leq p<\infty$ and 
$$\sum\limits_n|\hat S(\xi+n)|^p =\sum\limits_n|\hat {g_\xi^\#} 
(n)|^p\leq{\|\hat {h_\xi^\#}\|_{l_1}}^p{\|({S^\#})^\wedge\|_{l_p}}^p\leq 
C\|({S^\#})^\wedge\|_p^p$$
 where the constant $C$ does not depend upon $\xi$. So,
$\sup\limits_{\xi\in [0,1]}\sum\limits_n|\hat S(\xi+n)|^p < \infty$.

For $p=1$, by the Fourier inversion, we may assume  that $S$ is 
continuous. Now for a fixed $x$ define $g_x=e^{-2\pi ix.}S$. Then 
$g_x$ is continuous and $supp\;g_x\subseteq [\frac{1}{4},\frac{3}{4}]$.
 Also $\hat g_x^\#(n)=\hat S(x+n)$, where $g_x^\#$ is the 1-periodic 
 extension of $g_x$ from $[0,1)$. Therefore 
 $g_x^\#(t)=\sum\limits_{n\in\mathbb  Z}\hat S(x+n) e^{2\pi int}$ for 
 $t\in[0,1)$. As both sides of this equality are continuous functions 
 they will agree at $0$. So, 
 $g_x^\#(0)=\sum\limits_{n\in\mathbb  Z}\hat S(x+n)$ or 
 $S(0)=\sum\limits_{n\in\mathbb  Z}\hat S(x+n)$.

We will also need the following convolution result to prove our 
theorem.
\begin{theorem} Suppose $G$ is a locally compact abelian group. 
Let $1<r<2$. Then 
$L^r*L^{r^\prime}(\hat G)\subseteq M_p(\hat G)$
where  $\frac{2r}{3r-2}\leq p\leq \frac{2r}{2-r}$ and 
$\frac{1}{r}+\frac{1}{r^\prime}=1$.
\label{thm:three3}
\end{theorem}
The above result is given in ~\cite[page 126]{L}. 
The  main ingredient of the proof is the use of a Multilinear 
Riesz-Thorin  Interpolation theorem~\cite{zyg}.
\noindent {\bf Remark:}  For $r=2$, 
$L^2*L^2(G) \subseteq M_p(\hat G),\;\;\;\forall p\in [1,\infty)$.
\noindent{\bf Proof of Theorem~\ref{thm:two3}: }
Let $0<r<1$ and assume $1<p<2$. Define 
$k_r(x)=\sum\limits_{n\in\mathbb  Z}\phi(n)r^{|n|}e^{2\pi inx}$
 for $x\in [0,1)$. Then $k_r\in L^1(\mathbb  T)$ and 
 $\hat k_r(n)=\phi(n)r^{|n|}$. Thus 
$\hat k_r\in l_1(\mathbb  Z)$ and 
\begin{eqnarray*}
\|\hat k_r\|_{p^\prime}\leq \|\phi\|_{p^\prime}.
\end{eqnarray*}
Define $F_r(x)=k_r(x)S(x)$ for $x\in[0,1)$. Clearly 
$F_r\in L^1(\mathbb R)$ and 
$supp\;F_r\subseteq [\frac{1}{4},\frac{3}{4}]$  
 and $\hat F_r(\xi)=\sum\limits_n\phi(n)\hat S(\xi-n)r^{|n|}$. Then 
$\hat F_r|_\mathbb  Z(l)=\hat k_r*(S^\#)^\wedge(l)$, so 
by Theorem~\ref{thm:three3} $\hat F_r|_\mathbb  Z\in M_q(\mathbb  Z)$ 
for 
$q\in [\frac{2p}{3p-2},\frac{2p}{2-p}]$ with 
$\|\hat F_r|_\mathbb  Z\|_{M_q(\mathbb  Z)}\leq C_p\|\phi
\|_{p^\prime}\tau_p$. 
Hence,(by Theorem 3.4 of~\cite{BG}) $\hat F_r\in M_q(\mathbb R)$ for 
$q\in [\frac{2p}{3p-2},\frac{2p}{2-p}]$  with 
\begin{equation}
\|\hat F_r\|_{M_q(\mathbb R)}\leq C_p\tau_p
\|\phi\|_{p^\prime}
\label{eqn:two3}
\end{equation}
Again from Lemma~\ref{lem:one3} and dominated convergence theorem we 
have 
$\hat F_r(\xi)\rightarrow W_{\phi,\hat S}(\xi)\;{\mbox{a.e.}}\;\; 
{\mbox{as}}\;\;  r\rightarrow 1$.
 Therefore from the inequality~(\ref{eqn:two3}) we have 
 $W_{\phi,\hat S}\in M_q(\mathbb R)$
and $\|W_{\phi,\hat S}\|_{M_q(\mathbb R)}\leq C_p
\tau_p\|\phi\|_{p^\prime}$
 for $q\in [\frac{2p}{3p-2},\frac{2p}{2-p}]$. 
 Similarly for $2<p<\infty$ , by the same argument we get 
$$W_{\phi, \hat S}\in {M_q(\mathbb R)}$$ for 
$q\in [\frac{2p}{p+2},\frac{2p}{p-2}]$ and 
$$\|W_{\phi, \hat S}\|_{M_q(\mathbb R)}\leq 
C_p\tau_p\|\phi\|_{p^\prime}.$$  
This completes the proof of the theorem.

We will now relax the hypothesis that  
$supp\; S\subseteq[\frac{1}{4},\frac{3}{4}]$ to allow  $S$ to have 
arbitrary compact support by imposing a certain extra condition on $S$. 
For this we need the following lemma, which is easy to prove.
\begin{lem}
Let $A:\mathbb R\rightarrow \mathbb R$ defined by $A(x)=\alpha x$ for 
some $0\not =\alpha\in\mathbb  Z$. Then if 
$\Lambda\in S_p(\mathbb R)$ then $\Lambda\circ A\in S_p(\mathbb R)$.
\label{lem:add}
\end{lem}
Suppose $supp \;S\subseteq [-N,N]$ and  
$\sum\limits_{n\in\mathbb  Z}|\hat S(\xi+n)|^p<\infty$ for all 
$\xi\in[0,1)$. Define $S_N(x)=S(4Nx-2N)$. Then 
$supp\;S_N\subseteq[\frac{1}{4},\frac{3}{4}]$. Also from the 
condition on $\hat S$ we have 
$\sum\limits_{n\in\mathbb  Z}|(S_N^\#)^\wedge|^p<\infty$. Thus if 
$\phi\in l_{p^\prime}(\mathbb  Z)$ from Theorem~\ref{thm:two3} we have 
$W_{\phi,\hat S_N}\in M_q(\mathbb R)$ for the values of $q$ mentioned 
in the statement of the theorem. This along with Lemma~\ref{lem:add}  
says that  $W_{\phi,\hat S}\in M_q(\mathbb R)$. So in particular
\begin{cor}
Let $S\in C_C^1(\mathbb R)$ and $1<p<2$. Then for 
$\phi\in l_{p^\prime}$, $W_{\phi,\hat S}\in M_q(\mathbb R)$ for 
$q\in [\frac{2p}{3p-2},\frac{2p}{2-p}]$.
\end{cor}
By putting additional restrictions on $\phi$ we have the following 
(note that $l_p(\mathbb  Z)\subset l_{p^\prime}(\mathbb  Z)$)
\begin{prop}
Let $1<p<2$. Suppose $S\in L^p(\mathbb R)$ and has compact support. For 
$\phi\in l_p(\mathbb  Z)$ define 
$W_{\phi, \hat S}(\xi)=\sum\limits_n{\phi(n)\hat S(\xi-n)}$. Then 
$W_{\phi, \hat S}\in {M_q(\mathbb R)}$ for {\bf all} $1\leq q <\infty$
and $\|W_{\phi, \hat S}\|_{M_q(\mathbb R)}\leq C\|\phi\|_p\|S\|_p$.
\label{prop:one4}
\end{prop}
\noindent{\bf Proof:}
Let $supp\;S\subseteq [-N,N]$ for some $N\in \mathbb  N$. Define 
$S_N(x)=S(4Nx-2N)$. Then 
$supp\;S_N\subseteq [\frac{1}{4},\frac{3}{4}]$   and 
$S_N\in L^p(\mathbb R)$. Let 
$S_N^\#$ be 1-periodic extension of $S_N$ from $[0,1)$ and 
$(S_N^\#)^\wedge\in l_{p^\prime}(\mathbb  Z)$ (as $\hat S\in L^{p^\prime}(\mathbb R)$). Now,  if 
$\phi\in l_p(\mathbb  Z)$ then $\phi*(S_N^\#)^\wedge\in M_q(\mathbb  
Z)$ for $1\leq q<\infty$, 
because of the following reason. Consider the operator
\begin{eqnarray*}
T:l_1(\mathbb  Z)\times L^1(\mathbb  T)\longrightarrow M_q(\mathbb  Z) 
\\
T:l_2(\mathbb  Z)\times L^2(\mathbb  T)\longrightarrow M_q(\mathbb  Z),
\end{eqnarray*}
defined by $T(\phi,f)=\phi*\hat f$. Then 
$\|T(\phi,f)\|_{M_q(\mathbb  Z)}\leq\|\phi\|_{l_p
(\mathbb  Z)}\|f\|_{L^p(\mathbb  T)}$ for $p=1$ or 2. 
 So by Multilinear Riesz-Thorin interpolation  theorem~\cite{zyg}, 
 $T$ is a
 bounded and multilinear operator from 
 $l_p(\mathbb  Z)\times L^p(\mathbb  T)$ into 
$M_q(\mathbb  Z)$ for $1<p<2$ and $q\in[1,\infty)$. Thus 
$\phi*(S_N^\#)^\wedge\in M_q(\mathbb  Z)$. Following the same 
approach 
as in the proof of Theorem~\ref{thm:two3} we have 
$W_{\phi,\hat S_N}\in {M_q(\mathbb R)}$. So 
by Lemma~\ref{lem:add}  we have $W_{\phi,\hat S}\in {M_q(\mathbb R)}$ 
and
 $\|W_{\phi,\hat S}\|_{M_q(\mathbb R)}\leq C\|\phi\|_p\|S\|_p$,  where 
 $C$ is a constant
 depending on support of $S$.

\noindent{\bf{Remark:}}
Observe that our result (Theorem~\ref{thm:two3}) does   not 
match with Jodeit's result (Theorem~\ref{thm:one3}) in the 
limiting case $p=1$. In our case $\phi$ is just a bounded sequence 
but Jodeit considered $\phi$ to be in 
$M_q(\mathbb  Z)=l_\infty(\mathbb  Z)\cap M_q(\mathbb  Z).$ For 
the case $p=2$, we have $W_{\phi, \hat S}\in M_q(\mathbb R)$ for all 
$q\in[1,\infty)$ whenever $\phi\in l_2(\mathbb  Z)$. From Plancherel 
theorem it is easy to see that 
$l_2(\mathbb  Z)=l_2(\mathbb  Z)\cap M_q(\mathbb  Z)$ for all 
$q\in[1,\infty)$.  These observations pose the following problem:

`` Let $S\in L^1(\mathbb R)$, 
$supp\;S\subseteq [\frac{1}{4},\frac{3}{4}]$ and 
$\sum\limits_{n\in\mathbb  Z} |(S^\#)^\wedge(n)|^p<\infty$, for 
$1<p<\infty$. For $\phi\in l_p(\mathbb  Z)\cap M_q(\mathbb  Z)$ define 
$W_{\phi,\hat S}=\sum\limits_{n\in\mathbb  Z}\phi(n)\hat S(\xi-n)$. 
Then is it true that $W_{\phi,\hat S}\in M_q(\mathbb R)?"$

\end{document}